\documentclass[11pt]{article}
\usepackage{amssymb,amsfonts,amsmath,amsthm,cite}
\usepackage{epsfig}
\newtheorem{thm}{Theorem}[section]
\newtheorem{cor}{Corollary}[section]
\newtheorem{lem}{Lemma}[section]

\makeatletter \@addtoreset{equation}{section}

\def\pf{\noindent {\it Proof.\ }}
\def\qed{\hfill \rule{4pt}{7pt}}

\begin{document}

\begin{center}
{\large {\bf Congruences for Bipartitions with Odd Parts Distinct}}

\vskip 6mm

{\small William Y.C. Chen$^1$ and Bernard L.S. Lin$^2$
\\[%
2mm] Center for Combinatorics, LPMC-TJKLC\\
Nankai University, Tianjin 300071,
P.R. China \\[3mm]
$^1$chen@nankai.edu.cn, $^2$linlishuang@cfc.nankai.edu.cn \\[0pt%
] }
\end{center}

\noindent {\bf Abstract.}  Hirschhorn and Sellers studied arithmetic
properties of the number of partitions with odd parts distinct. In
another direction, Hammond and Lewis investigated arithmetic
properties of the number of bipartitions. In this paper, we consider
the number of bipartitions with odd parts distinct. Let this number
be denoted by  $pod_{-2}(n)$. We obtain two Ramanujan type
identities for $pod_{-2}(n)$, which imply that $pod_{-2}(2n+1)$ is
even  and $pod_{-2}(3n+2)$ is divisible by $3$. Furthermore, we show
that
 for any $\alpha\geq 1$ and
$n\geq 0$, $ pod_{-2}\left(3^{2\alpha+1}n+\frac{23\times
3^{2\alpha}-7}{8}\right) $ is a multiple of $3$ and
$pod_{-2}\left(5^{\alpha+1}n+\frac{11\times 5^\alpha+1}{4}\right)$
is divisible by $5$. We also find combinatorial interpretations for
the  two congruences modulo $2$ and $3$.


\noindent \textbf{Keywords:} partition, bipartition, congruence,
birank.

\noindent \textbf{AMS Classification:} 05A17, 11P83

\section{Introduction}
A partition $\lambda$ of a positive integer $n$ is any
non-increasing sequence of positive integers whose sum is $n$. The
weight of $\lambda$ is the sum of its parts, denoted by $|\lambda|$.
 A bipartition $\pi$
of $n$ is a pair of partitions $(\pi_1, \pi_2)$ with
$|\pi_1|+|\pi_2|=n$. Let $p_{-2}(n)$ denote the number of
bipartitions of $n$. The generating function for $p_{-2}(n)$ equals
\[
\sum_{n=0}^\infty p_{-2}(n)q^n=\frac{1}{(q;q)_\infty^2}.
 \]
In this paper, we shall employ the standard $q$-series notation 
\cite{Andrews98}
\[
(a;q)_0=1,\quad (a;q)_n:=\prod_{k=0}^{n-1}(1-aq^k), \quad \text{for
}n\geq 1,
\]
and
\[
(a_1,a_2,\ldots,a_m;q)_\infty=\lim_{n\to
\infty}\prod_{j=1}^m(a_j;q)_n,\quad |q|<1.
\]

The function $p_{-2}(n)$ has drawn much interest, see, for example, \cite{Atkin68, CJW06, Garvan09,
Gandhi63, HL04, Ramanathan50}.  Ramanathan \cite{Ramanathan50}
established the following congruences:
\begin{equation}\label{eqRamanathan}
p_{-2}(5n+2)\equiv p_{-2}(5n+3)\equiv p_{-2}(5n+4)\equiv 0\ ({\rm
mod\ }5),
\end{equation}
which are analogous to the classical congruences of Ramanujan, namely,
\begin{equation}\label{eqRamanujan5}
 p(5n+4)\equiv 0\ ({\rm mod\ }5)
\end{equation}
and
\begin{equation}\label{eqRamanujan7}
 p(7n+5)\equiv 0\ ({\rm mod\ }7),
\end{equation}
 where $p(n)$ is the number of partitions of $n$.

Dyson \cite{Dyson44} defined the rank of a partition as the largest
part minus the number of parts. Let $N(r,t,n)$ denote the number of
partitions of $n$ whose rank is congruent to $r$ modulo $t$.  Aktin
and Swinnerton-Dyer \cite{AS54}
 proved the following conjecture of Dyson \cite{Dyson44}
\[
N(r,5,5n+4)=\frac{p(5n+4)}{5}\quad 0\leq r\leq 4,
\]
and
\[
N(r,7,7n+5)=\frac{p(7n+5)}{7}\quad 0\leq r\leq 6.
\]

For a bipartition $\pi=(\pi_1,\pi_2)$,
Hanmmond and Lewis \cite{HL04} defined the birank $b(\pi)$ as
\begin{equation}\label{eqbirank1}
b(\pi)=n(\pi_1)-n(\pi_2),
\end{equation}
where $n(\lambda)$ denotes the number of parts of $\lambda$. It has
been shown that the birank $b(\pi)$ can be used to give
combinatorial interpretations of the congruences in
\eqref{eqRamanathan}. Recently, Garvan \cite{Garvan09} defined two
biranks. One can be utilized to explain all the three  congruences in \eqref{eqRamanathan}, while the other is valid for two of the three
congruences.

 We wish to consider bipartitions with odd parts
 distinct. Recall that Andrews, Hirschhorn and Sellers \cite{AHS09} have
  investigated  arithmetic properties of partitions with even parts
  distinct. Hirschhorn and Sellers \cite{HS10}
  considered arithmetic properties of partitions with odd parts
  distinct. To be precise, by a bipartition with odd parts distinct we mean a
bipartition $\pi=(\pi_1,\pi_2)$ for which  the odd parts of $\pi_1$
are distinct and the odd parts of $\pi_2$ are also distinct. Notice
that $\pi_1$ and $\pi_2$ are allowed to have an odd part in common.
For example, there are $11$ bipartitions of $4$:
\begin{eqnarray*}
& &((4),\emptyset)\ ((3,1),\emptyset) \ ((2,2),\emptyset)\
((3),(1))\ ((2,1),(1))\ ((2),(2))
\\[5pt]
& &((1),(2,1))\ ((1),(3))\ (\emptyset,(2,2))\ (\emptyset,(3,1))\
(\emptyset,(4)).
\end{eqnarray*}

Let $pod_{-2}(n)$ denote the number of  bipartitions of $n$ with odd
parts distinct. It is easy to derive the  generating function for $pod_{-2}(n)$, that is,
\begin{equation}\label{eqgenpod2}
\sum_{n=0}^\infty
pod_{-2}(n)q^n=\frac{(-q;q^2)_\infty^2}{(q^2;q^2)_\infty^2}.
\end{equation}

The main objective of this paper is to study  arithmetic properties
of $pod_{-2}(n)$ in the spirit of Ramanujan's congruences for the
partition function $p(n)$. We shall prove that
\begin{equation}
\sum_{n=0}^\infty pod_{-2}(2n+1)q^n=\frac{2(q^8;q^8)_\infty
^2}{(q;q)_\infty ^3(q^4;q^4)_\infty}
\end{equation}
and
\begin{equation}\label{eqgenpod3n+21}
\sum_{n=0}^\infty pod_{-2}(3n+2)q^n=3\frac{(q^2;q^2)_\infty
^4(q^6;q^6)_\infty ^6}{(q;q)_\infty ^6(q^4;q^4)_\infty^6},
\end{equation}
which implies that for all $n\geq 0$,
\begin{equation}\label{eqcongpod2n+1}
pod_{-2}(2n+1)\equiv 0\ ({\rm mod\ }2)
\end{equation}
and
\begin{equation}\label{eqcongpod3n+2}
pod_{-2}(3n+2)\equiv 0\ ({\rm mod\ }3).
\end{equation}

We also give three infinite families of congruences modulo $3$ and two
infinite families of congruences
 modulo $5$. For example,  for $\alpha\geq 1$ and $n\geq 0$,
\begin{equation}
pod_{-2}\left(3^{2\alpha+1}n+\frac{23\times
3^{2\alpha}-7}{8}\right)\equiv 0\ ({\rm mod\ }3)
\end{equation}
and
\begin{equation}
pod_{-2}\left(5^{\alpha+1}n+\frac{11\times 5^\alpha+1}{4}\right)
\equiv 0\ ({\rm mod\ }5).
\end{equation}

Furthermore, we  show that the birank $b(\pi)$ defined by Hammond
and Lewis can be used to explain the congruence
\eqref{eqcongpod3n+2}. Furthermore, we introduce another birank to
give a combinatorial explanation of \eqref{eqcongpod3n+2}. Our
birank $c(\pi)$ of a bipartition $\pi=(\pi_1,\pi_2)$ is defined by
\begin{equation}\label{eqbirank2}
c(\pi)=l(\pi_1)-l(\pi_2),
\end{equation}
where $l(\lambda)$ denotes the largest part of  $\lambda$. It is
worth mentioning that neither of the  two biranks $b(\pi)$ and
$c(\pi)$ leads to a combinatorial interpretation of the congruence
\eqref{eqcongpod2n+1}. It should  be noted  that the birank $c(\pi)$
is not exact the conjugate of $b(\pi)$ for bipartitions with odd
parts distinct because the conjugation of such a bipartition no
longer preserves this property.

This paper is organized as follows.  In Section 2, two  identities
of Ramanujan type are obtained. In Section 3, three infinite
families of congruences modulo $3$ for $pod_{-2}(n)$ are
established. In section 4, we obtain two  infinite families of congruences
modulo $5$ for $pod_{-2}(n)$. In Section 5, we prove
that both  biranks $b(\pi)$ and $c(\pi)$ can  be applied to give
a combinatorial interpretation
 of the fact that
$pod_{-2}(3n+2)$ is a multiple of $3$. We also give a simple
combinatorial  explanation of the fact that $pod_{-2}(2n+1)$ is even for any
$n$.

\section{Two Ramanujan-type identities}

In this section, we shall prove the following two
Ramanujan-type identities for the number of bipartitions with odd parts distinct.

\begin{thm}\label{thmgenofp2} We have
\begin{eqnarray}
\sum_{n=0}^\infty pod_{-2}(2n+1)q^n&=&\frac{2(q^8;q^8)_\infty
^2}{(q;q)_\infty ^3(q^4;q^4)_\infty},\label{eq2n+1ofp2}\\[5pt]
\sum_{n=0}^\infty pod_{-2}(3n+2)q^n&=&3\frac{(q^2;q^2)_\infty
^4(q^6;q^6)_\infty ^6}{(q;q)_\infty
^6(q^4;q^4)_\infty^6}.\label{eq3n+2ofp2}
\end{eqnarray}
\end{thm}

We need some properties of the function $\psi(q)$, namely,
 \begin{equation}\label{eqpsi}
 \psi(q)=\sum_{n=0}^\infty q^{n(n+1)/2}.
 \end{equation}
 Let
$f(a,b)$ be Ramanujan's general theta function given by
\[
f(a,b)=\sum_{n=-\infty}^\infty a^{n(n-1)/2}b^{n(n+1)/2},\quad
|ab|<1.
\]
 Jacobi's triple product
identity can be stated in Ramanujan's notation as follows
\begin{equation}\label{eqJacobi}
f(a,b)=(-a,ab)_\infty(-b;ab)_\infty(ab;ab)_\infty.
\end{equation}
 Thus,
\begin{equation}\label{eqpsiproduct}
\psi(-q)=f(-q,-q^3)=\frac{(q^2;q^2)_\infty}{(-q;q^2)_\infty}.
\end{equation}
 Combining \eqref{eqgenpod2} and \eqref{eqpsiproduct}, we
obtain that
\begin{equation}\label{eqgenpod2psi}
\sum_{n=0}^\infty pod_{-2}(n)q^n=\frac{1}{\psi(-q)^2}.
\end{equation}

It should be noted that Bringmann and Lovejoy \cite{BL06} have
studied  arithmetic properties of the numbers $\overline{pp}(n)$,
which are the coefficients $q^n$ in $1/\varphi(-q)^2$, namely,
\[
\sum_{n=0}^\infty \overline{pp}(n)q^n=\frac{1}{\varphi(-q)^2},
\]
where
 \begin{equation}\label{eqvarphi}
 \varphi(q)=\sum_{n=-\infty}^\infty q^{n^2}.
 \end{equation}

\begin{lem}\label{lem1} We have
\begin{eqnarray}
\frac{1}{\psi(-q)}&=&\frac{1}{(q^2;q^2)_\infty(q^4;q^4)_\infty}\left(f(q^6,q^{10})
+qf(q^2,q^{14})\right)\label{eqlem1of1}\\[5pt]
&=&\frac{\psi(-q^9)}{\psi(-q^3)^4}\left(A(-q^3)^2+
qA(-q^3)\psi(-q^9)+q^2\psi(-q^9)^2\right),\label{eqlem1of2}
\end{eqnarray}
where
\[
A(q)=\frac{(q^2;q^2)_\infty(q^3;q^3)_\infty ^2}{(q;q)_\infty
(q^6;q^6)_\infty}.
\]
\end{lem}
\pf It is easily checked that
\begin{equation}\label{eqlem1temp1}
\psi(q)\psi(-q)=\frac{(q^2;q^2)_\infty}{(q;q^2)_\infty}\cdot
\frac{(q^2;q^2)_\infty}{(-q;q^2)_\infty}=\frac{(q^2;q^2)_\infty^2}{(q^2;q^4)_\infty}
=(q^2;q^2)_\infty(q^4;q^4)_\infty.
\end{equation}
From \cite[Corollary (ii), p.49]{Berndt91}, it follows that
\begin{equation}\label{eqlem1temp2}
\psi(q)=f(q^6,q^{10})+qf(q^2,q^{14}).
\end{equation}
Dividing \eqref{eqlem1temp2} by \eqref{eqlem1temp1}, we are led to
the $2$-dissection \eqref{eqlem1of1} of $1/\psi(-q)$. The proof of
\eqref{eqlem1of2} is a little more involved; See \cite[Lemma
2.2]{HS10} for the details.\qed

In view of the above lemma, we are in a position to prove Theorem
\ref{thmgenofp2}.

 \noindent {\it Proof of Theorem \ref{thmgenofp2}.}
 By the
$2$-dissection \eqref{eqlem1of1} of $1/\psi(-q)$  and the generating
function \eqref{eqgenpod2psi} for $pod_{-2}(n)$, we see that
\begin{eqnarray*}
\sum_{n=0}^\infty pod_{-2}(n)q^n&=&
\frac{1}{(q^2;q^2)_\infty^2(q^4;q^4)_\infty^2}\left(f(q^6,q^{10})
+qf(q^2,q^{14})\right)^2.
\end{eqnarray*}
Considering the coefficients of $q^{2n+1}$ on both sides, we observe
that
\[
\sum_{n=0}^\infty
pod_{-2}(2n+1)q^n=\frac{2}{(q;q)_\infty^2(q^2;q^2)_\infty^2}f(q^3,q^5)f(q,q^7).
\]
Consequently, we get \eqref{eq2n+1ofp2},  since
\begin{eqnarray*}
f(q^3,q^5)f(q,q^7)&=&(-q,-q^3,-q^5,-q^7;q^8)_\infty(q^8;q^8)_\infty^2\\[5pt]
&=&(-q;q^2)_\infty(q^8;q^8)_\infty^2\\[5pt]
&=&\frac{(q^2;q^2)_\infty^2(q^8;q^8)_\infty^2}{(q;q)_\infty(q^4;q^4)_\infty}.
\end{eqnarray*}
This completes the proof of \eqref{eq2n+1ofp2}.

 By the $3$-dissection \eqref{eqlem1of2} of
$1/\psi(-q)$, we find that
\[
\sum_{n=0}^\infty
pod_{-2}(n)q^n=\frac{\psi(-q^9)^2}{\psi(-q^3)^8}\left(A(-q^3)^2+
qA(-q^3)\psi(-q^9)+q^2\psi(-q^9)^2\right)^2.
\]
Extracting the terms $q^{3n+2}$ on both sides, we obtain
\[
\sum_{n=0}^\infty
pod_{-2}(3n+2)q^{3n+2}=3q^2\frac{\psi(-q^9)^2}{\psi(-q^3)^8}A(-q^3)^2\psi(-q^9)^2.
\]
By dividing both sides $q^2$ and replacing
$q^3$ by $q$, arrive at
\eqref{eq3n+2ofp2}. This completes the proof. \qed

As  consequences of Theorem
\ref{thmgenofp2}, we obtain the following congruences.

\begin{cor}\label{cor1}For each nonnegative integer $n$,
\[
pod_{-2}(2n+1)\equiv 0\ ({\rm mod\ }2)\quad \text{and} \quad
pod_{-2}(3n+2)\equiv 0\ ({\rm mod\ }3).
\]
\end{cor}

\section{Three infinite families of congruences modulo $3$ }

In this section, we  wish to establish the following three infinite
families of Ramanujan-like congruences modulo $3$ satisfied by
$pod_{-2}(n)$ by two different approaches. The proof of Theorem
\ref{thmmod3ofp2} needs the formula for  the number of ways to represent an
integer $n$
as a sum of two triangular numbers as well as
 a characterization of integers that can not be written as a sum of
two squares.  On the other hand,  Theorem \ref{thmmod32ofp2}  follows from
the generating function for the numbers $ pod_{-2}(3n+1)$.
For notational convenience, we assume that all the congruences in this
section are modulo $3$.

\begin{thm}\label{thmmod3ofp2}
 For all $\alpha\geq 0$ and $n\geq 0$,
\begin{equation}\label{eq3^n1}
pod_{-2}\left(3^{2\alpha+1}n+\frac{23\times
3^{2\alpha}-7}{8}\right)\equiv 0\ ({\rm mod\ }3).
\end{equation}
\end{thm}

\begin{thm}\label{thmmod32ofp2}
 For all $\alpha\geq 1$ and $n\geq 0$,
\begin{equation}\label{eq3^n2}
pod_{-2}\left(3^{2\alpha+1}n+\frac{7\times
3^{2\alpha}+1}{4}\right)\equiv 0\ ({\rm mod\ }3)
\end{equation}
and
\begin{equation}\label{eq3^n3}
pod_{-2}\left(3^{2\alpha+1}n+\frac{11\times
3^{2\alpha}+1}{4}\right)\equiv 0\ ({\rm mod\ }3).
\end{equation}
\end{thm}

To prove the above congruences,   the following
lemma is useful.

\begin{lem}\label{lem2}
\begin{eqnarray}
\psi(q)&=&f(q^3,q^6)+q\psi(q^9),\label{eqlem2of1}\\[5pt]
\psi(q^3)&\equiv& \psi(q)^3.\label{eqlem2of2}
\end{eqnarray}
\end{lem}
\pf From \cite[Corollary (ii), p.49]{Berndt91} it is clear that the
identity \eqref{eqlem2of1} holds. Since
\[
(1-q^n)^3\equiv (1-q^{3n})\ ({\rm mod\ }3)
\]
and
\[
\psi(q)=\frac{(q^2;q^2)_\infty^2}{(q;q)_\infty}=\frac{\prod_{n\geq
1}(1-q^{2n})^2}{\prod_{n\geq 1}(1-q^n)},
\]
we obtain  \eqref{eqlem2of2}. This completes the proof.  \qed

 \noindent {\it Proof of Theorem \ref{thmmod3ofp2}.}  By Lemma
\ref{lem2}, we have
\[
\sum_{n=0}^\infty pod_{-2}(n)(-q)^n=\frac{\psi(q)}{\psi(q)^3}\equiv
\frac{f(q^3,q^6)+q\psi(q^9)}{\psi(q^3)}.
\]
Extracting the terms $q, q^4, q^7,\ldots$ on both sides of above identity,
dividing by $q$, and replacing $q^3$ by $q$, we get
\begin{equation}\label{eqiden3n+1}
\sum_{n=0}^\infty (-1)^{n+1}pod_{-2}(3n+1)q^n\equiv
\frac{\psi(q^3)}{\psi(q)} \equiv \psi(q)^2.
\end{equation}
Let the numbers  $t_2(n)$ be defined by
\[
\psi(q)^2=\sum_{n=0}^\infty t_2(n)q^n.
\]
 By comparing the coefficients
of $q^n$ on both sides of \eqref{eqiden3n+1}, we find that for each
$n\geq 0$,
\begin{equation}\label{eqcongt2n}
pod_{-2}(3n+1)\equiv (-1)^{n+1}t_2(n).
\end{equation}
From \cite[Theorem 3.6.2]{Berndt06}, it follows that for all  integers
$n\geq 0$,
\[
t_2(n)=d_{1,4}(4n+1)-d_{3,4}(4n+1),
\]
where $d_{j,k}(n)$ denotes the number of positive divisors $d$ of
$n$ such that $d\equiv j\ ({\rm mod\ }k)$. Moreover, by
\cite[Theorem 2.15]{MNZ91}, we have that $d_{1,4}(n)-d_{3,4}(n)=0$
if and only if every prime $p$ congruent to $3$ modulo $4$ in the
canonical factorization of $n$ appears with an odd exponent.

 It is clear that  for $\alpha\geq 1$ and $n\geq 0$, the integer $s=4\times
3^{2\alpha}n+\frac{23\times 3^{2\alpha-1}-3}{2}$ is a multiple of $3$
but not divisible by $9$.  This implies that
\[
t_2\left(\frac{s-1}{4}\right)=d_{1,4}(s)-d_{3,4}(s)=0.
\]
Substituting $n=\frac{s-1}{4}$ into \eqref{eqcongt2n}, we obtain that
\[
pod_{-2}\left(3^{2\alpha+1}n+\frac{23\times
3^{2\alpha}-7}{8}\right)\equiv 0\ ({\rm mod\ }3).
\]
The case $\alpha=0$ has been considered in Corollary \ref{cor1}. This completes the
proof.\qed

 \noindent {\it Proof of Theorem \ref{thmmod32ofp2}.}
Invoking the identity \eqref{eqiden3n+1} in the proof of Theorem
\ref{thmmod3ofp2}, we deduce that
\begin{equation}\label{eqiden3n+1sec}
\sum_{n=0}^\infty pod_{-2}(3n+1)q^n\equiv -\psi(-q)^2.
\end{equation}
Applying \eqref{eqlem1of2} and \eqref{eqlem2of2} to
\eqref{eqiden3n+1sec}, we obtain
\begin{eqnarray*}
\sum_{n=0}^\infty
pod_{-2}(3n+1)q^n&=&-\frac{\psi(-q)^3}{\psi(-q)}\equiv
-\frac{\psi(-q^3)}{\psi(-q)}\\[5pt]
&\equiv& -\frac{\psi(-q^9)}{\psi(-q^3)^3}\left(A(-q^3)^2+
qA(-q^3)\psi(-q^9)+q^2\psi(-q^9)^2\right).\\[5pt]
\end{eqnarray*}
Extracting the  terms $q^{3n+2}$ for $n\geq 0$, we find that
\[
\sum_{n=0}^\infty pod_{-2}(9n+7)q^{3n+2}\equiv
-q^2\frac{\psi(-q^9)^3}{\psi(-q^3)^3}.
\]
Dividing both sides of above identity by $q^2$ and replacing $q^3$
by $q$, we see that
\[
\sum_{n=0}^\infty pod_{-2}(9n+7)q^{n}\equiv
-\frac{\psi(-q^3)^3}{\psi(-q)^3}\equiv -\psi(-q^3)^2.
\]
Similarly, it can be shown that
\begin{equation}\label{eqiden27n+7}
\sum_{n=0}^\infty pod_{-2}(27n+7)q^n\equiv -\psi(-q)^2
\end{equation}
and for $n\geq 0$,
\[
pod_{-2}(27n+16)\equiv pod_{-2}(27n+25)\equiv 0.
\]
So the proof is complete for the case $\alpha=1$.
Combining \eqref{eqiden3n+1sec} and \eqref{eqiden27n+7}, it can be seen that
 for $n\geq 0$,
\begin{equation}\label{eqcongpodrelation}
pod_{-2}(3n+1)\equiv pod_{-2}(27n+7).
\end{equation}
By induction on $\alpha$, it is easy to establish congruences
\eqref{eq3^n2} and \eqref{eq3^n3} based on the relation
\eqref{eqcongpodrelation}. \qed

\section{Two infinite families of congruences modulo $5$}

In this section, we  give two infinite
families of Ramanujan-like congruences modulo $5$ satisfied by
$pod_{-2}(n)$ from a modular equation of degree $5$ due to
Ramanujan. For notational convenience, we assume that all the
congruences in this section are modulo $5$.

\begin{thm}\label{thmmod5ofp2}
 For all $\alpha\geq 1$ and $n\geq 0$,
\begin{equation}\label{eq5^n1}
pod_{-2}\left(5^{\alpha+1}n+\frac{11\times 5^\alpha+1}{4}\right)
\equiv 0\ ({\rm mod\ }5)
\end{equation}
and
\begin{equation}\label{eq5^n1}
pod_{-2}\left(5^{\alpha+1}n+\frac{19\times 5^\alpha+1}{4}\right)
\equiv 0\ ({\rm mod\ }5).
\end{equation}
\end{thm}

To prove the above congruences, we need the following lemma.
\begin{lem}\label{lem3}
Let $1\leq r\leq 4$. Let the numbers $a(n)$ be given by
\[
\sum_{n=0}^\infty a(n)q^n=\sum_{n=0}^\infty
\frac{q^{5n+r}}{1-q^{10n+2r}}.
\]
Then \[\sum_{n=0}^\infty a(5n)q^n=\sum_{n=0}^\infty
\frac{q^{5n+r}}{1-q^{10n+2r}}.
\]
\end{lem}
\pf Clearly,
\begin{equation}\label{eqmod5temp1}
\sum_{n=0}^\infty \frac{q^{5n+r}}{1-q^{10n+2r}}=\sum_{n=0}^\infty
\sum_{k=0}^\infty q^{(5n+r)(2k+1)}.
\end{equation}
Since for $1\leq r\leq 4$ and $k\geq 0$, $(5n+r)(2k+1)$ is a
multiple of $5$ if and only if $k\equiv 2\ ({\rm mod\ }5)$. It
follows that
\begin{eqnarray*}
\sum_{n=0}^\infty a(5n)q^{5n}&=&\sum_{n=0}^\infty \sum_{k\equiv 2\
({\rm mod\ }5)}^\infty q^{(5n+r)(2k+1)}\\[5pt]
&=&\sum_{n=0}^\infty \sum_{t=0}^\infty q^{(5n+r)(10t+5)}.
\end{eqnarray*}
Replacing $q^5$ by $q$ and using \eqref{eqmod5temp1}, we complete the proof. \qed

\noindent {\it Proof of Theorem \ref{thmmod5ofp2}.} It is easy to
deduce the following relation
\begin{equation}\label{eqmod5temp2}
\psi(q^5)\equiv \psi(q)^5.
\end{equation}
From the generating function \eqref{eqgenpod2psi} for  $pod_{-2}(n)$
and \eqref{eqmod5temp2} it follows that
\begin{eqnarray}
q\sum_{n=0}^\infty pod_{-2}(n)(-q)^n&=&\frac{q}{\psi(q)^2}
\equiv \frac{q\psi(q)^3\psi(q^5)}{\psi(q^5)^2}\nonumber\\[5pt]
&\equiv&
\frac{q\psi(q)^3\psi(q^5)-5q^2\psi(q)\psi(q^5)^3}{\psi(q^5)^2}.\label{eqmod5temp3}
\end{eqnarray}
From  \cite[Entry 8(i), p.249]{Berndt91}, we see that
\begin{eqnarray*}
q\psi(q)^3\psi(q^5)-5q^2\psi(q)\psi(q^5)^3&=&\sum_{n=0}^\infty
\frac{(5n+1)q^{5n+1}}{1-q^{10n+2}}-\sum_{n=0}^\infty
\frac{(5n+2)q^{5n+2}}{1-q^{10n+4}}\\[5pt]
& & \;\;-\sum_{n=0}^\infty \frac{(5n+3)q^{5n+3}}{1-q^{10n+6}}+
\sum_{n=0}^\infty \frac{(5n+4)q^{5n+4}}{1-q^{10n+8}}.
\end{eqnarray*}
This implies that
\begin{eqnarray*}
q\psi(q)^3\psi(q^5)-5q^2\psi(q)\psi(q^5)^3&\equiv&\sum_{n=0}^\infty
\frac{q^{5n+1}}{1-q^{10n+2}}-\sum_{n=0}^\infty
\frac{2q^{5n+2}}{1-q^{10n+4}}\\[5pt]
& & \;\;-\sum_{n=0}^\infty \frac{3q^{5n+3}}{1-q^{10n+6}}+
\sum_{n=0}^\infty \frac{4q^{5n+4}}{1-q^{10n+8}}.
\end{eqnarray*}
Write the above series modulo $5$ as
\[ \sum_{n=0}^\infty A(n)q^n.\]
Applying Lemma \ref{lem3} yields
\[
\sum_{n=0}^\infty A(5n)q^n\equiv
q\psi(q)^3\psi(q^5)-5q^2\psi(q)\psi(q^5)^3.
\]
Extracting the  terms $q^{5n}$  from
\eqref{eqmod5temp3} and replacing $q^5$ by $q$, we have
\[
-\sum_{n=0}^\infty
pod_{-2}(5n+4)(-q)^{n+1}\equiv\frac{\sum\limits_{n=0}^\infty
A(5n)q^n}{\psi(q)^2}.
\]
Combining the above two equations, we find that
\begin{eqnarray}
-\sum_{n=1}^\infty pod_{-2}(5n-1)(-q)^{n}&\equiv&
\frac{q\psi(q)^3\psi(q^5)-5q^2\psi(q)\psi(q^5)^3}{\psi(q)^2}\nonumber\\[5pt]
&\equiv&\frac{q\psi(q)^3\psi(q^5)}{\psi(q)^2}\nonumber\\[5pt]
&\equiv&q\psi(q)\psi(q^5)\label{eqmod5temp4}\\[5pt]
&=&\psi(q^5)\left(qf(q^{10},q^{15})+q^2f(q^5,q^{20})+q^4\psi(q^{25})\right).
\label{eqmod5temp5}
\end{eqnarray}
Note that the last equation follows from \cite[Corollary (ii),
p.49]{Berndt91}. Comparing coefficients of $q^{5n+a}(a=0,3,4)$
in \eqref{eqmod5temp5}, we see that for $n\geq 0$,
\begin{equation}\label{eqmod5temp6}
pod_{-2}(25n+14)\equiv pod_{-2}(25n+24)\equiv 0
\end{equation}
and
\[
\sum_{n=0}^\infty pod_{-2}(25n+19)(-q)^{n+1}\equiv
q\psi(q)\psi(q^5).
\]
In view of the above identity and \eqref{eqmod5temp4}, we deduce that
\[
\sum_{n=0}^\infty pod_{-2}(25n+19)(-q)^{n+1}\equiv
-\sum_{n=0}^\infty pod_{-2}(5n+4)(-q)^{n+1},
\]
which implies that for  $n\geq 0$,
\[
pod_{-2}(25n+19)\equiv -pod_{-2}(5n+4).
\]
  Using the above relation and \eqref{eqmod5temp6}, it is easily checked by
  induction that for $\alpha\geq 1$,
\[
pod_{-2}\left(5^{\alpha+1}n+\frac{11\times 5^\alpha+1}{4}\right)
\equiv pod_{-2}\left(5^{\alpha+1}n+\frac{19\times
5^\alpha+1}{4}\right)\equiv 0.
\]
This completes the proof.\qed

It should be noted that Chan \cite{Chan08} has used modular forms to
establish  infinite families of congruences modulo $5$ for
Andrews-Paule's $2$-diamond partitions.  His approach can also be used
to prove the congruences in this section.

\section{Combinatorial interpretations }

In this section, we show that both the  biranks $b(\pi)$ and
$c(\pi)$ can be used to give a combinatorial interpretation of the
fact that $pod_{-2}(3n+2)$ is divisible by $3$. We conclude this paper with
a simple explanation of the parity of $pod_{-2}(2n+1)$.

Let $R(m,n)$ denote the number of bipartitions $\pi$ of $n$ with odd
parts distinct such that birank $b(\pi)=m$. By using the
transformation that interchanges $\pi_1$ and $\pi_2$ in
\eqref{eqbirank1}, we see that
\begin{equation}
R(m,n)=R(-m,n).
\end{equation}
Let $R(r,t,n)$ be the number of bipartitions $\pi$ of $n$ with odd
parts distinct such that birank $b(\pi)$ is congruent to $r$ modulo $t$,
i.e.,
\[
R(r,t,n)=\sum_{m\equiv r\ ({\rm mod\ }t)} R(m,n).
 \]
 Then we have $R(r,t,n)=R(t-r,t,n)$. Moreover, it is easy to derive the following
 generating function for $R(m,n)$,
\begin{equation}\label{eqgenbirank1}
\sum_{n=0}^\infty \sum_{m=-\infty}^\infty
R(m,n)z^mq^n=\frac{(-qz;q^2)_\infty
(-q/z;q^2)_\infty}{(q^2z;q^2)_\infty(q^2/z;q^2)_\infty}.
\end{equation}
The above formula enables us to obtain generating function for the numbers of the
form $R(r,t, n)-R(s,t,n)$.

\begin{thm} \label{thmgenofbirank}
\begin{eqnarray}
 \sum_{n=0}^\infty
\left(R(0,2,n)-R(1,2,n)\right)q^n&=&\frac{\varphi(-q)}{\psi(q^2)},\label{eqgenbirank2n}\\[5pt]
\sum_{n=0}^\infty
\left(R(0,3,n)-R(1,3,n)\right)q^n&=&\frac{\psi(-q)}{\psi(-q^3)}\label{eqgenbirank3n},\\[5pt]
\sum_{n=0}\left(R(0,4,n)-R(2,4,n)\right)q^n&=&
\frac{\varphi(q^2)}{\psi(q^2)}.\label{eqgenbirank4n}
\end{eqnarray}
\end{thm}

\pf Taking $z=-1$ in the generating function \eqref{eqgenbirank1},
we have
\begin{eqnarray*}
\sum_{n=0}^\infty
\left(R(0,2,n)-R(1,2,n)\right)q^n&=&\frac{(q;q^2)_\infty^2}{(-q^2;q^2)_\infty
^2}=\frac{(q;q)_\infty^2}{(q^4;q^4)_\infty^2}\\[5pt]
&=&\frac{(q;q)_\infty^2}{(q^2;q^2)_\infty}\times
\frac{(q^2;q^4)_\infty}{(q^4;q^4)_\infty}\\[5pt]
&=&\frac{\varphi(-q)}{\psi(q^2)}.
\end{eqnarray*}
Note that the last equation holds since
\[
\varphi(q)=f(q,q)=(-q;q^2)_\infty^2(q^2;q^2)_\infty.
\]
 Substituting $z=\xi=e^{2\pi i/3}$ into both sides of
\eqref{eqgenbirank1} and applying the relation $R(1,3,n)=R(2,3,n)$, we
find that
\begin{eqnarray*}
\frac{(-q\xi;q^2)_\infty
(-q\xi^2;q^2)_\infty}{(q^2\xi;q^2)_\infty(q^2\xi^2;q^2)_\infty}&=&
\sum_{n=0}^\infty
\left(R(0,3,n)+R(1,3,n)\xi+R(2,3,n)\xi^2\right)q^n\\[5pt]
&=&\sum_{n=0}^\infty \left(R(0,3,n)-R(1,3,n)\right)q^n.\\[5pt]
\end{eqnarray*}
Since  $1-x^3=(1-x)(1-x\xi)(1-x\xi^2)$, we see that
\[
\frac{(-q\xi;q^2)_\infty
(-q\xi^2;q^2)_\infty}{(q^2\xi;q^2)_\infty(q^2\xi^2;q^2)_\infty}=
\frac{(-q^3;q^6)_\infty(q^2;q^2)_\infty}{(q^6;q^6)_\infty(-q;q^2)_\infty}
=\frac{\psi(-q)}{\psi(-q^3)}.
\]
Hence we arrive at the relation
\eqref{eqgenbirank3n}. Similarly,  setting  $z=i$ in
\eqref{eqgenbirank1} and using the fact that $R(1,4,n)=R(3,4,n)$, we
get
\begin{eqnarray*}
\sum_{n=0}^\infty\left(R(0,4,n)-R(2,4,n)\right)q^n=\frac{(-qi,qi;q^2)_\infty}
{(q^2i,-q^2i;q^2)_\infty}
=\frac{(-q^2;q^4)_\infty}{(-q^4;q^4)_\infty}.
\end{eqnarray*}
 It remains to show that
\begin{eqnarray*}
\frac{\varphi(q^2)}{\psi(q^2)}&=&(-q^2;q^4)_\infty^2(q^4;q^4)_\infty\times
\frac{(q^2;q^4)_\infty}{(q^4;q^4)_\infty}\\[5pt]
&=&(-q^2;q^4)_\infty(q^4;q^8)_\infty=\frac{(-q^2;q^4)_\infty}{(-q^4;q^4)_\infty}.
\end{eqnarray*}
This completes the proof.\qed

Based on the relation \eqref{eqgenbirank3n}, we see that
 the birank given by Hammond and
Lewis leads to a classification of the bipartitions  of $3n+2$ with odd parts distinct
into three equinumerous sets. Thus we deduce the following theorem.

\begin{thm}\label{thmbirank}
For  $0\leq r \leq 2$ ,
\[
R(r,3,3n+2)=\frac{pod_{-2}(3n+2)}{3}.
\]
\end{thm}

\pf By \eqref{eqlem2of1} and \eqref{eqgenbirank3n}, we find
\[
\sum_{n=0}^\infty
\left(R(0,3,n)-R(1,3,n)\right)q^n=\frac{f(-q^3,q^6)-q\psi(-q^9)}{\psi(-q^3)}.
\]
Since the term $q^{3n+2}$ does not appear on
the right-hand side of above identity, it follows that for $n\geq 0$,
\[
R(0,3,3n+2)=R(1,3,3n+2).
\]
Combining the fact that $R(1,3,n)=R(2,3,n)$ and the relation
$\sum_{r=0}^2R(r,3,n)=pod_{-2}(n)$, we conclude that for $0\leq
r\leq 2$,
\[
R(r,3,3n+2)=\frac{pod_{-2}(3n+2)}{3}.
\]
This completes the proof.\qed

We now use the new birank $c(\pi)$ to give another interpretation of
the congruence relation \eqref{eqcongpod3n+2}. As above, we need to
consider $R_2(m,n)$ as the number of bipartitions $\pi$ of $n$ with
odd parts  distinct and whose birank $c(\pi)$ equals  $m$. The
following theorem gives the generating function for $R_2(m,n)$.

\begin{thm}
\begin{equation}\label{eqgenbirank2}
\sum_{n=0}^\infty \sum_{m=-\infty}^\infty R_2(m,n)z^mq^n=
(1+q/z)(1+qz)\frac{(-q^3z^2;q^2)_\infty}{(q^2z^2;q^2)_\infty}
\frac{(-q^3/z^2;q^2)_\infty}{(q^2/z^2;q^2)_\infty}.
\end{equation}
\end{thm}
\pf Let $A_k(n)$ (resp. $B_k(n)$) denote the number of partitions of
$n$ such that the  odd parts are distinct and  the largest part
equals $2k$ (resp. $2k+1$). It is easy to see that
\[
A(z,q):=\sum_{k=0}^\infty \sum_{n=0}^\infty A_k(n)z^{2k}q^n
=\sum_{m=0}^\infty\frac{q^{2m}z^{2m}(-q;q^2)_m}{(q^2;q^2)_m}
\]
and
\[
B(z,q):=\sum_{k=0}^\infty \sum_{n=0}^\infty B_k(n)z^{2k+1}q^n
=\sum_{m=0}^\infty\frac{q^{2m+1}z^{2m+1}(-q;q^2)_m}{(q^2;q^2)_m}.
\]
By the $q$-binomial theorem \cite[Theorem 1.3.1]{Berndt06}, we find
that\[ A(z,q)=\frac{(-q^3z^2;q^2)_\infty}{(q^2z^2;q^2)_\infty}
\]
and
\[
B(z,q)=qz\frac{(-q^3z^2;q^2)_\infty}{(q^2z^2;q^2)_\infty}.
\]
 Let
$\pi=(\lambda,\mu)$ be a bipartition. Consider the parities of the largest parts
of $\lambda$ and $\mu$. We have
\begin{eqnarray*}
\sum_{n=0}^\infty \sum_{m=-\infty}^\infty R_2(m,n)z^mq^n&=&
A(z,q)A(1/z,q)+A(z,q)B(1/z,q)\\[5pt]
& & \;\;+B(z,q)A(1/z,q)+B(z,q)B(1/z,q)\\[5pt]
&=&(1+q/z+qz+q^2)\frac{(-q^3z^2;q^2)_\infty}{(q^2z^2;q^2)_\infty}
\frac{(-q^3/z^2;q^2)_\infty}{(q^2/z^2;q^2)_\infty}\\[5pt]
&=&(1+q/z)(1+qz)\frac{(-q^3z^2;q^2)_\infty}{(q^2z^2;q^2)_\infty}
\frac{(-q^3/z^2;q^2)_\infty}{(q^2/z^2;q^2)_\infty}.
\end{eqnarray*}
This completes the proof. \qed

 Let $R_2(r,t,n)$
denote the number of bipartitions $\pi$ of $n$  with  odd parts
distinct such that  $c(\pi)\equiv r\ ({\rm mod\ }t)$. We are now
ready to show that the  birank $c(\pi)$ implies a  combinatorial
explanation of congruence \eqref{eqcongpod3n+2}.

\begin{thm}For $0\leq r \leq 2$,
\[
R_2(r,3,3n+2)=\frac{pod_{-2}(3n+2)}{3}.
\]
\end{thm}

\pf Let $\xi=e^{2\pi i/3}$. Substituting $z=\xi$ into
\eqref{eqgenbirank2}, we get
\begin{eqnarray*}
\sum_{n=0}^\infty \sum_{r=0}^2 R(r,3,n)\xi^r q^n
&=&(1+q\xi)(1+q\xi^2)\frac{(-q^3\xi^2;q^2)_\infty}{(q^2\xi^2;q^2)_\infty}
\frac{(-q^3\xi;q^2)_\infty}{(q^2\xi;q^2)_\infty}\\[5pt]
&=&\frac{(-q\xi^2;q^2)_\infty}{(q^2\xi^2;q^2)_\infty}
\frac{(-q\xi;q^2)_\infty}{(q^2\xi;q^2)_\infty}\\[5pt]
&=&\frac{(-q^3;q^6)_\infty}{(q^6;q^6)_\infty}\sum_{n=0}^\infty
(-q)^{n(n+1)/2}.
\end{eqnarray*}
Since no triangular numbers
$n(n+1)/2$ are congruent to $2$ modulo $3$, equating coefficients of $q^n$ on both sides, we find that, for each
integer $n\geq 0$,
\[
\sum_{r=0}^2 R(r,3,3n+2)\xi^r =0.
 \]
 Consequently,
 \[
 R_2(0,3,3n+2)=R_2(1,3,3n+2)=R_2(2,3,3n+2),
 \]
 since $\xi$ is one of the roots of the irreducible polynomial
 $1+z+z^2=0$. This  completes the proof.\qed

 Here is a simple combinatorial  explanation of the congruence
 \eqref{eqcongpod2n+1}. Let us  define the rank $d(\pi)$ of a bipartition
 $\pi=(\pi_1,\pi_2)$ as
 the number of parts of $\pi_1$.
  Let $R_3(m,n)$ denote the number of bipartitions
 $\pi$ of $n$ with odd parts distinct and whose rank $d(\pi)=m$. Let $R_3(r,t,n)$
denote the number of bipartitions $\pi$ of $n$  with  odd parts
distinct such that  $d(\pi)\equiv r\ ({\rm mod\ }t)$. The
 generating function for $R_3(m,n)$ equals
 \[
 \sum_{n=0}^\infty \sum_{m=0}^\infty R_3(m,n)z^mq^n=
 \frac{(-qz;q^2)_\infty}{(q^2z;q^2)_\infty}\cdot
 \frac{(-q;q^2)_\infty}{(q^2;q^2)_\infty}.
 \]
 Setting $z=-1$, we get
 \begin{equation}\label{eqrank2n}
 \sum_{n=0}^\infty \left(R_3(0,2,n)-R_3(1,2,n)\right)q^n=
 \frac{(q^2;q^4)_\infty}{(q^4;q^4)_\infty},
 \end{equation}
 which immediately implies that
 \begin{equation}\label{eqrank2n1}
 R_3(0,2,2n+1)=R_3(1,2,2n+1)
 \end{equation}
 and
 \begin{equation}\label{eqrank2n2}
\sum_{n=0}^\infty \left(R_3(0,2,2n)-R_3(1,2,2n)\right)(-q)^n=
 \frac{(-q;q^2)_\infty}{(q^2;q^2)_\infty}.
 \end{equation}
In light of \eqref{eqrank2n1},  we see that the rank $d(\pi)$ leads to a combinatorial interpretation of the
congruence \eqref{eqcongpod2n+1}.
 It is worth mentioning  that the right-hand side of  \eqref{eqrank2n2}
 is the generating function for partitions with odd parts distinct
 \cite{HS10}.

\vspace{.2cm} \noindent{\bf Acknowledgments.} This work was
supported by the 973 Project, the PCSIRT Project of the Ministry of
Education, and the National Science Foundation of China.

\end{document}